\documentclass[10pt]{amsart}

\usepackage[utf8]{inputenc}
\usepackage[T1]{fontenc}

\usepackage{times,amsfonts,amsmath,amstext,amsbsy,amssymb,
  amsopn,amsthm,upref,eucal, amscd}
\usepackage{colonequals}
\usepackage{hyperref}
\usepackage{cleveref}

\usepackage[pdftex]{graphicx}

\usepackage[style=alphabetic,natbib=true,backend=biber,giveninits=true]{biblatex}

\newtheorem{theorem}{Theorem}[section]

\numberwithin{equation}{section}

\theoremstyle{definition}

\newtheorem{remark}[theorem]{Remark}
\newtheorem{conjecture}[theorem]{Conjecture}

\def\leq{\leqslant }
\def\geq{\geqslant}

\DeclareFieldFormat{title}{\mkbibquote{#1}}

\begin{document}

\title{Lower bounds on the dimension of the Rauzy gasket}

\author{Rodolfo Gutiérrez-Romo}
\address{Rodolfo Gutiérrez-Romo: Institut de Mathématiques de Jussieu -- Paris Rive Gauche, UMR 7586, Bâtiment Sophie Germain, 75205 Paris Cedex 13, France.}
\email{rodolfo.gutierrez@imj-prg.fr}
\urladdr{http://rodol.fo}

\author{Carlos Matheus}
\address{Carlos Matheus: Centre de Mathématiques Laurent Schwartz, CNRS (UMR 7640), École Polytechnique, 91128 Palaiseau, France.}
\email{carlos.matheus@math.cnrs.fr}
\urladdr{http://carlos.matheus.perso.math.cnrs.fr}


\begin{abstract}
The Rauzy gasket $R$ is the maximal invariant set of a certain renormalization procedure for special systems of isometries naturally appearing in the context of Novikov's problem in conductivity theory for monocrystals. 

It was conjectured by Novikov and Maltsev in 2003 that the Hausdorff dimension $\dim_{\mathrm{H}}(R)$ of Rauzy gasket is strictly comprised between $1$ and $2$. 

In 2016, Avila, Hubert and Skripchenko confirmed that $\dim_{\mathrm{H}}(R)<2$. In this note, we use some results by Cao--Pesin--Zhao in order to show that $\dim_{\mathrm{H}}(R)>1.19$. 
\end{abstract}
\maketitle


\section{Introduction}

\begin{figure}[t!]
	\includegraphics[width=0.5\textwidth]{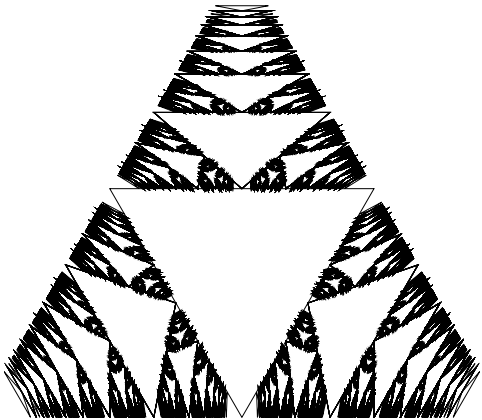}
	\caption{The Rauzy gasket.}
	\label{fig:rauzy_gasket}
\end{figure}

The Rauzy gasket is a fractal subset of the standard $2$-simplex related to frequencies of letters in ternary episturmian words \cite{AS}, dynamics of special systems of isometries \cite{D}, and a particular case of Novikov's problem around the trajectories of electrons on Fermi surfaces in the presence of constant magnetic fields \cite{dLD,AHS2}. The Rauzy gasket is depicted in \Cref{fig:rauzy_gasket}.
 
Concretely, the Rauzy gasket is defined as follows. Consider the standard $2$-simplex $\Delta = \{(x_1, x_2, x_3)\in\mathbb{R}^3_+: x_1+x_2+x_3 = 1\}$. We decompose $\Delta$ into three simplices $\Delta_j = \{(x_1,x_2,x_3)\in\Delta: x_j\geq\sum_{k\neq j} x_k\}$ and a hole $\Delta\setminus\bigcup_{j=1}^3\Delta_j$. The projectivizations of the matrices
\[M_1 = \left(\begin{array}{ccc} 1 & 1 & 1 \\ 0 & 1 & 0 \\ 0 & 0 & 1\end{array}\right), \quad M_2 = \left(\begin{array}{ccc}1&0&0 \\ 1&1&1 \\ 0&0&1\end{array}\right), \quad M_3 = \left(\begin{array}{ccc}1&0&0 \\ 0&1&0 \\ 1&1&1\end{array}\right)\]
induce weakly contracting maps $f_j\colon\Delta\to\Delta_j$, $j=1,2,3$.

In this context, recall from \cite{AS} that the \emph{Rauzy gasket} is the unique non-empty compact subset of $\Delta$ such that
\[R = f_1(R)\cup f_2(R) \cup f_3(R).\]

The fact that the Rauzy gasket has zero Lebesgue measure was proved by several authors including Levitt \cite{L}\footnote{Using an argument attributed to Yoccoz.}, Arnoux--Starosta \cite{AS} and de Leo--Dynnikov \cite{dLD}. 

It was conjectured by Novikov and Maltsev \cite{NM} in 2003 that:

\begin{conjecture}[Novikov--Maltsev] $1<\dim_{\mathrm{H}}(R)<2$.
\end{conjecture}

Some numerical experiments by R. de Leo and I. Dynnikov \cite{dLD} suggest that $1.7<\dim_{\mathrm{H}}(R)<1.8$, and Avila--Hubert--Skripchenko \cite{AHS} established that $\dim_{\mathrm{H}}(R)<2$.

The main result of this note is the following theorem:

\begin{theorem}\label{t.A} $\dim_{\mathrm{H}}(R) > 1.19$.
\end{theorem}

The proof of this result occupies the remainder of this text. 

\section{Lower bounds on the Hausdorff dimension of the Rauzy gasket}

In this section, we give a lower bound on $\dim_{\mathrm{H}}(R)$ via the construction of appropriate uniformly expanding repellers inside $R$. 

\subsection{General framework}

We will use somewhat general methods to obtain bounds for the Hausdorff dimension of a uniformly expanding repeller in dimension $2$. These methods rely on estimating the singular values of the derivatives of the maps defining the set. More precisely, given $n$ uniformly contracting maps $T_1, \dotsc, T_n\colon X \to X$, where $X \subseteq \mathbb{R}^2$ is a compact set, and a repeller $K$ defined as the unique non-empty compact set such that $K = \bigcup_{k = 1}^n T_k(K)$, we need to estimate quantities of the form $\max_{x \in X}\|D_x T_k\|$ and $\min_{x \in X}\|(D_x T_k)^{-1}\|$, where $\|\cdot\|$ denotes the largest singular value. Since for any $a, b, c, d \in \mathbb{R}$ one has
\[
	\left\| \begin{pmatrix}
		a & b \\
		c & d
	\end{pmatrix} \right\| = \sqrt{\frac{a^2 + b^2 + c^2 + d^2 + \sqrt{(a^2 + b^2 + c^2 + d^2)^2 - 4(ad - bc)^2}}{2}},
\]
we obtain the simple estimates
\[
	\sqrt{\frac{a^2 + b^2 + c^2 + d^2}{2}} \leq \left\| \begin{pmatrix}
		a & b \\
		c & d
	\end{pmatrix}\right\| \leq \sqrt{a^2 + b^2 + c^2 + d^2},
\]
which we will write as $\|\cdot\|^- \leq \|\cdot\| \leq \|\cdot\|^+$.


There are several methods in the literature to obtain lower bounds on the Hausdorff dimension of repellers. For our purposes, the thermodynamical method of Cao--Pesin--Zhao \cite{CPZ} is quite useful. In a nutshell, they consider a repeller $\Lambda$ of a $C^2$-expanding map $g$ on a surface, a parameter $1\leq s\leq 2$, and the potential $\psi^s(x,g) = \log\alpha_1(x,g)+(s-1)\log\alpha_2(x,g)$, where $\alpha_1(x,g)\geq \alpha_2(x,g)$ are the singular values of $D_x g$. Observe that $\|D_x g\|^- \leq \alpha_1(x, g) \leq \|D_x g\|^+$ and $\|(D_x g)^{-1}\|^- \leq \alpha_2(x, g)^{-1} \leq \|(D_x g)^{-1}\|^+$.

By Corollary 3.1 of \cite{CPZ}, one has that 
\[\dim(\Lambda)\geq s_1\]
where $s_1$ is the unique root of the equation $P(g,-\psi^s(\cdot ,g))=0$ and $P(g,\theta)$ stands for the topological pressure of the potential $\theta$, i.e.,  
\[P(g,\theta)\colonequals\sup\left\{h_{\mu}(g)+\int\theta\, d\mu(x): \mu \text{ is } g\text{-invariant} \right\}\]
(see (3.2) and (2.4) in \cite{CPZ}). The theory of (subadditive) thermodynamical formalism (as explained\footnote{Cf.\ Lemma 3.2 of \cite{Fa} in particular.} in Section 3 of \cite{Fa}, for instance) states that 
\[P(g,\theta)<0 \iff \sum_{m\geq 1}\sum_{x\in \textrm{Fix}(g^m)} \exp(\theta_m(x))<\infty\]
where $\theta_m(x)\colonequals\sum_{j=0}^{m-1}\theta(g^j(x))$. 

In general, $s\mapsto P(g, -\psi^s(\cdot,g))$ is a continuous and strictly decreasing function of $s$. Therefore, $s_1\geq s_0$ for all $s_0$ with \[\sum\limits_{m\geq 1}\sum_{x\in \textrm{Fix}(g^m)} \exp\left(-\sum_{j=0}^{m-1}\psi^{s_0}(g^j(x),g)\right) = \infty.\]

\subsection{The Rauzy gasket}

Observe that each composition $f_k\circ f_j$, with $k\neq j$, is a contraction on $\Delta$ (cf.\ Lemma 2 in \cite{AS}). Thus, for each integer $n\geq 2$, the unique non-empty compact subset $K_n$ such that 
\[K_n = \bigcup\limits_{\substack{i \in S_n}} f_{i_n}\circ\dotsb\circ f_{i_1}(K_n),\]
where $S_n = \{1,2,3\}^n \setminus \{(1, \dotsc, 1), (2, \dotsc, 2), (3, \dotsc, 3)\}$, is a uniformly expanding repeller contained in $R$. 

In what follows, we consider the Riemannian metric on $T\Delta = \{(v_1,v_2,v_3)\in\mathbb{R}^3: v_1+v_2+v_3=0\}$ induced by the usual Euclidean scalar product of $\mathbb{R}^3$ normalized so that the vectors $(\varepsilon_1,\varepsilon_2,\varepsilon_3)$, $\{\varepsilon_1,\varepsilon_2,\varepsilon_3\} = \{-1,0,1\}$ have norm $1$. In particular, $\mathcal{B} = \{(1,-1,0), (-1,-1,2)/\sqrt{3}\}$ is an orthonormal basis of $T\Delta$. 

\begin{remark}
  A natural alternative is to consider the Fubini--Study metric $d(\mathbb{R}x,\mathbb{R}y) = \frac{\|x\wedge y\|}{\|x\| \|y\|}$ on the projective space $P\mathbb{R}^3$. However, we chose the \emph{ad hoc} Riemannian metric above because the operation of taking exterior powers would lead to heavier calculations.  
\end{remark}

The repeller $K_{13}$ defined by $g_{13}$ sending each $\Delta_i = f_{i_{13}}\circ\dotsb\circ f_{i_1}(\Delta)$, with $i \in S_{13}$, onto $\Delta$ is uniformly expanding with respect to this Riemannian metric. Indeed, we can estimate the smallest expansion factor as $1/\max_{i \in S_{13}} \max_{x \in \Delta_i} \|(D_x g_{13})^{-1}\|^+$ to obtain a value of at least $\sqrt{3}$.

Now, denote by
\[
	a=\log\left(\max_{i \in S_{13}}\max_{x \in \Delta_i}\|D_x g_{13}\|^+\right), \quad b=\log\left(1\mathbin{\Big/}\min_{i \in S_{13}}\min_{x \in \Delta_i}\|(D_x g_{13})^{-1}\|^-\right)
\]
and $\textrm{Fix}(g_{13}^m)=\exp(cm)$ for all $m$ (i.e., $c=\log(|S_{13}|) = \log(3^{13}-3)$). Observe that
\[
	\log \alpha_1(x, g) \leq a \quad \text{ and } \quad \log \alpha_2(x, g) \leq b
\]
for every $x \in \bigcup_{i \in S_{13}} \Delta_i$. Hence, 
\[\sum_{j=0}^{m-1}\psi^{s_0}(g_{13}^j(x),g_{13})\leq (a+b(s_0-1))m\]
and we deduce that  
\[\sum_{m\geq 1}\sum_{x\in \textrm{Fix}(g_{13}^m)} \exp\left(-\sum_{j=0}^{m-1}\psi^{s_0}(g_{13}^j(x),g_{13})\right)\geq \sum_{m\geq 1}\exp((c-a-b(s_0-1))m) = \infty\]
if $c-a-b(s_0-1)>0$, i.e., $s_0<1+(c-a)/b$.

In this way, obtain the bound
\[\dim_{\mathrm{H}}(K_{13})\geq s_1\geq 1+\frac{c-a}{b}.\]
With the help of a computer, we can find the exact values of $a$ and $b$. We obtain:
\[
	a = \log\left( 3208 \sqrt{ \frac{86185}{3} } \right), \quad b = \log\left( 4917248\sqrt{\frac{2}{1595}} \right) \quad \text{ and } \quad c = \log(1594320),
\]
which yields $\dim_{\mathrm{H}}(K_{13}) \geq 1 + \frac{c - a}{b} > 1.08$.

This lower bound can be improved by restricting to a smaller fractal. Indeed, instead of using every sequence in $S_{13}$, we can take a subset of such sequences designed to optimize the previous bound by decreasing the values of $a$ and $b$ while trying to maintain a large value of $c$. The heuristic we use is as follows:
\begin{enumerate}
	\item Sort the $i \in S_{13}$ according to $\max_{x \in \Delta_i} \|D_x g_{13}\|^+$ in an ascending order, assigning a number $r^+(i)$ to each $i \in S_n$.
	\item Find the $i \in S_{13}$ that maximizes $\log(r^+(i)) - \log(\max_{x \in \Delta_i} \|D_x g_{13}\|^+)$ and denote it by $i^*$. Let $S_{13}^+ = \{ i \in S_n : r^+(i) \leq r^+(i^*) \}$. For the remaining steps, we ignore the elements of $S_{13} \setminus S_{13}^+$. Let $a' = \log \max_{x \in \Delta_{i^*}} \|D_x g_{13}\|^+$.
	\item Sort the $i \in S_{13}^+$ according to $\min_{x \in \Delta_i} \|(D_x g_{13})^{-1}\|^-$ in a descending order, assigning a number $r^-(i)$ to each $i \in S_{13}^+$.
	\item Find the $i \in S_{13}^+$ that maximizes $\frac{\log(r^-(i)) - a'}{\log(1/\min_{x \in \Delta_i} \|(D_x g_{13})^{-1}\|^-)}$ and denote it by $i^{**}$. Our new set of sequences is now $S_{13}^{+-} = \{ i \in S_{13}^+ : r^-(i) \leq r^-(i^{**}) \}$ and we have that $b = \log(1/\min_{x \in \Delta_{i^{**}}} \|(D_x g_{13})^{-1}\|^-)$ and $c = \log(r^-(i^{**}))$. We also define $a = \log \max_{i \in S_{13}^{+-}} \max_{x \in \Delta_i} \|D_x g_{13}\|^+$ (it may happen that $a < a'$ as we have removed more sequences).
\end{enumerate}
We repeat this heuristic until the list of sequences does not change. We get the following final values:
\[
	a = \log\left( 6800 \sqrt{ \frac{829}{3} } \right), \quad b = \log\left( 615627\sqrt{\frac{3}{515}} \right) \quad \text{ and } \quad c = \log(898224).
\]
Thus, we obtain the bound $\dim_{\mathrm{H}}(K_{13}) \geq 1 + \frac{c - a}{b} > 1.19$ and we establish the lower bound $\dim_{\mathrm{H}}(R) > 1.19$.


\sloppy\printbibliography


\end{document}